\documentclass[11pt,english]{amsart}
\usepackage[T1]{fontenc}
\usepackage[utf8]{inputenc}
\usepackage{mathrsfs, amsmath, amstext, amsthm, amssymb, mathtools, physics}
\usepackage[all]{xy}
\usepackage{babel}
\usepackage{caption}
\makeatletter

\numberwithin{equation}{section}
\numberwithin{figure}{section}

\theoremstyle{plain}
\newtheorem{thm}{\protect\theoremname}[section]
\theoremstyle{remark}
\newtheorem{rem}[thm]{\protect\remarkname}
\theoremstyle{plain}

\theoremstyle{definition}

\theoremstyle{remark}

\theoremstyle{plain}
\newtheorem{lem}[thm]{\protect\lemmaname}
\theoremstyle{plain}
\newtheorem{cor}[thm]{\protect\corollaryname}
\theoremstyle{remark}

\newtheoremstyle{notestyle}
 {\topsep}
 {\topsep}
 {\normalfont}
 {}
 {\bfseries}
 {.}
 {.5em}
 {}
\theoremstyle{notestyle}

\usepackage{graphicx}
\usepackage{wrapfig}
\usepackage{fancyhdr}
\usepackage{extarrows}
\usepackage{geometry}
\makeatother

\providecommand{\definitionname}{Definition}
\providecommand{\examplename}{Example}
\providecommand{\lemmaname}{Lemma}
\providecommand{\corollaryname}{Corollary}
\providecommand{\propositionname}{Proposition}
\providecommand{\remarkname}{Remark}
\providecommand{\theoremname}{Theorem}

\global\long\def\CC{\mathbb{C}}%
\global\long\def\RR{\mathbb{R}}%

\global\long\def\G{\mathcal{G}}%
\global\long\def\rr{\boldsymbol{\underline{r}}}%
\DeclarePairedDelimiter\pInternal{(}{)}

\newcommand{\pr}[1]{\pInternal*{#1}}

\global\long\def\F{\mathcal{F}}%
\global\long\def\E{\mathcal{E}}%
\global\long\def\N{\mathcal{N}}%
\global\long\def\W{\mathcal{W}}%
\global\long\def\inv{^{\ast}}%
\global\long\def\su{\subseteq}%
\global\long\def\gl{\lambda}%

\global\long\def\g{\mathfrak{g}}%
\global\long\def\WW{\mathfrak{W}}%

\DeclareMathOperator{\Sym}{Sym}

\DeclareMathOperator{\Sing}{Sing}
\global\long\def\Da{\mathcal{D}_{\mathrm{arch}}}%
\DeclareMathOperator{\rk}{rk}

\newcommand{\KL}{\operatorname{KL}}
\newcommand{\RLCT}{\operatorname{RLCT}}
\newcommand{\Id}{\operatorname{Id}}

\usepackage{hyperref}

\title{A note on the learning coefficient of exceptional singularities}
\author{Giacomo Graziani}
\address{Department of Civil, Environmental and Architectural Engineering, University of Padua, 
Via F. Marzolo 9, 
35131 Padova, Italy.}
\email{giacomo.graziani@unipd.it}

\begin{document}

\begin{abstract}
We compute the learning coefficient of a polynomial neural network at a point of its neurovariety all of whose real parameters are good, and show that it equals $d/2$, with multiplicity one, where $d$ is the difference between the dimension of the parameter space and the dimension of the group of symmetries of the architecture. This is the value attained at a generic point of the neurovariety; consequently, exceptional singularities are invisible to the pair $(\lambda,\theta)$ whenever all real parameters in their fibre are good. By contrast, at a full parameter in the interior of the parameter domain whose differential has non-maximal rank we show that $\lambda<d/2$. The two statements together give a statistical reading of the geometric distinction between nodes and cusps.
\end{abstract}
\maketitle


\section*{Introduction}
 
To an algebraic neural network $X$ one associates its realisation map
\[
\Phi_X:\W\longrightarrow V,
\]
sending a choice of weights to the polynomial map it computes, and the neurovariety $\F(X)=\overline{\Phi_X(\W)}$. In \cite{Gra26b} we showed that neurovarieties carry singular points of a kind not produced by any degeneration of the architecture: these are the points of
\[
\E(X)=\pr{\Sing\F(X)\cap\Phi_X(\W)}\setminus\Da(X),
\]
called exceptional singularities, where $\Da(X)$ is the closure of the image of the non-full parameters. At such a point the neurovariety is singular, although it lies outside $\Da(X)$, so that every parameter realising it is full. The singularity may arise either from ramification of the realisation map or from the meeting of distinct branches, rather than from a degeneration of the network architecture.
 
On the statistical side, singular learning theory attaches to a statistical model a rational number $\gl>0$ and an integer $\theta\ge1$, the learning coefficient and its multiplicity, which govern the asymptotics of the free energy \eqref{eq:free-energy},
\[
F_n=nL_n(w_0)+\gl\log n-(\theta-1)\log\log n+O_p(1).
\]
For a regular model $\gl$ is half the number of parameters; the more singular the model, the smaller $\gl$, and the fewer effective degrees of freedom the model spends. The invariant $\gl$ is thus the natural bridge between the geometry of the fibres of $\Phi_X$ and the behaviour of learning, and one expects the learning coefficient to drop at the singular points of $\F(X)$.
 
This expectation is false, and the way in which it fails is the content of this note. Our main result is the following; see Theorem \ref{thm:A} and Corollary \ref{cor:dichotomy} for the precise statements, and \eqref{eq:dims} for the definition of $d$.
 
\noindent\textbf{Theorem.}
\emph{Let $u_0\notin\Da(X)$ be realised by a full real parameter in the interior of the parameter domain and assume that every real parameter of $\Phi_X^{-1}(u_0)$ in the parameter domain is good. Then}
\[
  \gl=\frac d2,\qquad \theta=1,
  \qquad
  d=\dim_\RR\W_\RR-\dim_\RR\G(X)(\RR).
\]
Since an exceptional singularity lies outside $\Da(X)$ and is realised by full parameters, the theorem applies to it whenever its fibre consists of good parameters, and the conclusion is then identical to the one at a generic point of $\F(X)$, that is the singularity leaves no trace on $(\gl,\theta)$. 
 
That this gives a genuine dichotomy is the content of Theorem \ref{thm:B}: at a full parameter in the interior of the parameter domain whose differential drops rank to $r<d$ one has
\[\gl\le\frac{d+r}{4}<\frac d2,\]
and this can be seen as a statistical reading of the geometric trichotomy of \cite[Theorem 2.7]{Gra26b}, that is, cusps are visible and nodes are not.

The key observation of the paper is that, at a good parameter, the germ of $\Phi_X$ is analytically equivalent to a linear projection, so the fibre ideal is generated by part of a coordinate system. No invariant computed from the germ at a single parameter can therefore distinguish such a point from one lying over a generic point of $\F(X)$. 

\section{Preliminaries}
\subsection{Thresholds and learning coefficients}
\label{Sec:rlct}

\subsubsection{Real log canonical thresholds of ideals}

Let $\Omega\su\RR^{m}$ be a compact semianalytic set with non-empty interior, let $f_1,\ldots,f_\ell$ be real analytic on a neighbourhood of $\Omega$, let $I=(f_1,\ldots,f_\ell)$ and let $\varphi\ge0$ be nearly analytic on a neighbourhood of $\Omega$, that is $\varphi=\varphi_a\varphi_s$ with $\varphi_a$ real analytic and $\varphi_s$ smooth and positive \cite[Section 1]{Lin17}. If the common zero locus $V_\Omega(I)$ of the $f_i$ in $\Omega$ is non-empty, the zeta function
\begin{equation}\label{eq:zeta}
\zeta\left(s\right)=\int_{\Omega}\pr{f_1^2+\cdots+f_\ell^2}^{-s/2}\varphi\left(w\right) \dd w
\end{equation}
converges for real $s$ smaller than some $\gl_1>0$ and extends to a meromorphic function on $\CC$ whose poles are positive rational numbers \cite[Corollary 2.6]{Lin17}. Following \cite[Section 1]{Lin17} the real log canonical threshold $\RLCT_\Omega(I;\varphi)=(\gl_1,\theta_1)$ is the smallest pole $\gl_1$ of $\zeta$ together with its order $\theta_1$; it depends only on $I$ and not on the generators \cite[Proposition 3.6]{Lin17}. Pairs are ordered as in \cite{Lin17}:
\[
(\gl_1,\theta_1)\preceq(\gl_2,\theta_2)
\quad\text{if}\quad \gl_1<\gl_2,\ \text{ or }\ \gl_1=\gl_2\ \text{and}\ \theta_1\ge\theta_2 ,
\]
and every minimum of pairs is taken with respect to $\preceq$. For $w\in\Omega$ we write $\RLCT_{\Omega,w}(I;\varphi)$ for the threshold over $\Omega\cap N$, where $N$ is any sufficiently small neighbourhood of $w$; it does not depend on $N$ \cite[Lemma 2.4]{Lin17} and the remark after \cite[Proposition 2.5]{Lin17}. If $w\in\operatorname{int}\Omega$ we simply write $\RLCT_{w}(I;\varphi)$.

We shall use the following properties.
\begin{enumerate}
  \item[(P1)] (\cite[Proposition 3.2(c)]{Lin17}). $$\RLCT_\Omega(I;\varphi)=\min_{w\in V_\Omega(I)}\RLCT_{\Omega,w}(I;\varphi).$$
  \item[(P2)] (\cite[Lemma 2.4]{Lin17}). If $\varphi\left(w\right)>0$ then $$\RLCT_{\Omega,w}(I;\varphi)=\RLCT_{\Omega,w}(I;1).$$ 
  In particular local thresholds at interior points are invariant under real analytic changes of coordinates, the Jacobian being absorbed in the amplitude (\cite[Proposition 3.8]{Lin17}).
  \item[(P3)] (\cite[Corollary 3.5]{Lin17}). If $c\sum_i f_i^2\le\sum_j h_j^2\le C\sum_i f_i^2$ on a neighbourhood of $w$ in $\Omega$ for some $c,C>0$, then 
  \[
  \RLCT_{\Omega,w}\left( \left(f_1,\dots,f_n\right);\varphi\right) = \RLCT_{\Omega,w}\left( \left(h_1,\dots,h_k\right);\varphi\right).
  \] 
  \item[(P4)] (\cite[Proposition 3.3]{Lin17}). If $w\in\partial\Omega$ and $W$ is a sufficiently small compact neighbourhood of $w$ in $\RR^m$, then $$\RLCT_{W}(I;\varphi)\preceq\RLCT_{\Omega,w}(I;\varphi).$$
  \item[(P5)] (\cite[Corollary 4.15]{Lin17}). If $w\in\operatorname{int}\Omega$, $w\in V_\Omega(I)$ and the Jacobian matrix of $(f_1,\ldots,f_\ell)$ at $w$ has rank $r$, then $$\RLCT_w(I;1)\preceq\pr{\frac{m+r}{2},1};$$ in particular the first entry of $\RLCT_w(I;1)$ is at most $(m+r)/2$.
\end{enumerate}

\subsubsection{The learning coefficient}
\label{sec:Learning}

Let $p(y\mid w)$ be a statistical model with parameter $w$ in a compact set $\Theta$, let $\varphi$ be a prior on $\Theta$ and let $Y_1,\ldots,Y_n$ be independent samples from $p(\,\bullet\mid w_0)$ for some $w_0\in\Theta$. Under the hypotheses of \cite[Chapter 6]{Wat09}, we have for the free energy
\[
F_n=-\log\int_{\Theta}\prod_{i=1}^{n}p(Y_i\mid w)\,\varphi\left(w\right)\dd w
\]
the expansion
\begin{equation}\label{eq:free-energy}
F_n =  L_n(w_0)\cdot n +\gl\cdot\log n - (\theta-1)\cdot \log\log n +O_p(1)
\end{equation}
where $L_n(w_0)=-\frac1n\sum_{i=1}^n\log p(Y_i\mid w_0)$ is the empirical
entropy at the true parameter and $(\gl,\theta)$, the learning coefficient and its multiplicity, are the smallest pole and its order, respectively, of $$z\mapsto\int_\Theta K\left(w\right)^{-z}\varphi\left(w\right)\dd w ,$$ with $K$ the Kullback-Leibler divergence from the true distribution \cite[Chapter 6]{Wat09}.
For a regular model, that is if the true parameter $w_0$ is unique, lies in the interior of $\Theta$, $\nabla^2K(w_0)\succ0$ and $\varphi(w_0)>0$, one has $\gl=\frac12\dim\Theta$ and $\theta=1$ \cite[Theorem 4]{Wat18}, \cite[Corollary 4.14]{Lin17}.

\begin{rem}\label{rem:factor-two}
If $K=c\sum_i f_i^2$ with $c>0$, then $\int_\Theta K^{-z}\varphi=c^{-z}\zeta(2z)$, with $\zeta$ as in \eqref{eq:zeta}, so that
\[
(2\gl,\theta)=\RLCT_\Theta\pr{(f_1,\ldots,f_\ell);\varphi}.
\]
This is the form taken by \cite[Theorem 1.2]{Lin17} when the Kullback-Leibler divergence is exactly a sum of squares, which is the case of the model of Section \ref{Sec:model}; in general \cite[Proposition 4.4]{Lin11} gives the same conclusion whenever $K$ is the composition of the parametrisation with a function with positive definite Hessian at the true parameter. Accordingly, for $w\in V_\Theta(K)$ we call half the first entry of $\RLCT_{\Theta,w}\pr{(f_1,\ldots,f_\ell);\varphi}$ the local learning coefficient at $w$, cf. \cite{LFWMW23}.
\end{rem}

\subsection{Algebraic networks as statistical models}
\label{Sec:model}

\subsubsection{Algebraic networks}
We recall what is needed from \cite{Gra26b}. Let $X=(\N,\W,\rr)$ be an algebraic network \cite[Definition 1.1]{Gra26b}, with realisation map
\[
\Phi_X:\W\longrightarrow V=X_{\mathrm{out}}\otimes\Sym^R\pr{X_{\mathrm{in}}\inv},\qquad R=r_1\cdots r_{L-1},
\]
neurovariety $\F(X)=\overline{\Phi_X(\W)}$, locus of full parameters $\W^\circ$ \cite[Definition 1.5]{Gra26b} and architectural degeneracy locus $\Da(X)=\overline{\Phi_X\pr{\W\setminus\W^\circ}}$. If $u_0\notin\Da(X)$ then $\Phi_X^{-1}(u_0)\su\W^\circ$. The group of symmetries $\G(X)$ \cite[(1.4)]{Gra26b} acts freely on $\W^\circ$ \cite[Remark 1.15]{Gra26b}, $\Phi_X$ is $\G(X)$-invariant, and $\Phi_X=\phi_X\circ\pi_X$ factors through the quotient $\pi_X:\W\to\WW(X)$ \cite[Theorem 1.12]{Gra26b}.

To keep the statements simple we assume from now on that $X$ is fully connected and $\W^\circ\neq\emptyset$; then $X$ satisfies \cite[(1.9)]{Gra26b} in view of \cite[Corollary 1.17]{Gra26b}, and $\W_\RR=\RR^{\dim\W}$. Everything below holds verbatim for networks satisfying \cite[(1.9)]{Gra26b}, as long as the parameter domain $\Theta$ lies in the smooth locus of $\W_\RR$. We set
\begin{equation}\label{eq:dims}
g=\dim\G(X),\qquad d=\dim\W-g=\dim\WW(X),\qquad k=\dim\F(X),
\end{equation}
the third equality by \cite[Theorem 1.16]{Gra26b}. Since the rank of $\dd \Phi_{X,w}$ is lower semicontinuous and $\F(X)$ is the closure of the image,
\begin{equation}\label{eq:max-rank}
  \rk \dd \Phi_{X,w}\le k\quad\text{for every }w\in\W,
\end{equation}
with equality on a dense open subset; for $w\in\W_\RR$ the rank is the same over $\RR$ and over $\CC$. For $w\in\W^\circ$ one has $\g(X)\cdot w\su\ker \dd \Phi_{X,w}$ with $\dim\g(X)\cdot w=g$, so $\rk \dd \Phi_{X,w}\le d$, and $w$ is good if equality holds \cite[Proposition 1.22 and Definition 1.23]{Gra26b}. If a good parameter exists then $k=d$, that is $X$ is generically finitely identifiable \cite[Corollary 1.24]{Gra26b}.

The exceptional locus is \cite[Definition 2.1]{Gra26b} $$\E(X)=\pr{\Sing\F(X)\cap\Phi_X(\W)}\setminus\Da(X).$$ 
A point $u\in\F(X)$ is a cusp if some $w\in\Phi_X^{-1}\left(u\right)$ has $\rk \dd \Phi_{X,w}<k$, a node if $\phi_X^{-1}\left(u\right)$ is not connected, and we write $\mathcal P(X)$ for the set of points around which $\phi_X$ is not proper \cite[Definition 2.4]{Gra26b}. In view of \cite[Theorem 2.7]{Gra26b} every exceptional singularity is a cusp, a node, or lies in $\mathcal P(X)$. If $X$ is generically finitely identifiable and $u\notin\Da(X)$, then $u$ is a cusp if and only if some parameter over $u$ is not good \cite[Section 2.1]{Gra26b}.

Since only real parameters enter the statistical model, we say that $u_0$ is a real cusp in $\Theta$ if some $w\in\Phi_X^{-1}(u_0)\cap\Theta$ has $\rk \dd \Phi_{X,w}<k$.

\subsubsection{The Gaussian model}

An admissible parameter domain is a compact set
$\Theta=\{g_1\ge0,\ldots,g_l\ge0\}\su\W_\RR$, with $g_j$ real polynomials, which is the closure of its interior. An amplitude is a smooth function $\varphi$ on a neighbourhood of $\Theta$ which is strictly positive on $\Theta$. Gaussian priors and the uniform prior on $\Theta$ are amplitudes. Fix $\sigma>0$ and an inner product on $V_\RR$, and consider the observation of the realised map in Gaussian noise,
\begin{equation}\label{eq:network-density}
p(y\mid w)=q\pr{y\mid \Phi_X\left(w\right)}, \qquad q(y\mid u)=\pr{2\pi\sigma^2}^{-M/2}\exp\pr{-\frac{\|y-u\|^2}{2\sigma^2}}, \qquad M=\dim_\RR V_\RR ,
\end{equation}
with true parameter $w_0\in\Theta$ and $u_0=\Phi_X(w_0)$. The base model $q$ is a Gaussian family, hence regular with Fisher matrix $\sigma^{-2}\Id_M$ and a direct computation gives
\begin{equation}\label{eq:network-KL}
K_{u_0}\left(w\right)=\KL\pr{q(\,\bullet\mid u_0)\,\Vert\,q(\,\bullet\mid \Phi_X\left(w\right))}=\frac{1}{2\sigma^2}\bigl\|\Phi_X\left(w\right)-u_0\bigr\|^2 .
\end{equation}
Hence, writing $I_{u_0}=\pr{\Phi_{X,1}-u_{0,1},\ldots,\Phi_{X,M}-u_{0,M}}$ in coordinates on $V_\RR$, Remark \ref{rem:factor-two} and Appendix \ref{App:hypotheses} give
\begin{equation}\label{eq:lambda-fibre}
(2\gl,\theta)=\RLCT_\Theta\pr{I_{u_0};\varphi} \overset{(P1)}{=} \min_{w\in\Phi_X^{-1}(u_0)\cap\Theta}\RLCT_{\Theta,w}\pr{I_{u_0};\varphi}.
\end{equation}
The ideal $I_{u_0}$ is the fibre ideal of \cite[Definition 1.8]{Lin11} and it does not depend on the chosen inner product.

\section{Main results}
\label{Sec:main}

\begin{thm}\label{thm:A}
Let $\Theta$ be an admissible parameter domain, let $u_0\in\Phi_X(\operatorname{int}\Theta)$ and let $\varphi$ be an amplitude. If $u_0$ is not a real cusp in $\Theta$, that is if $\rk \dd \Phi_{X,w}=k$ for every $w\in\Phi_X^{-1}(u_0)\cap\Theta$, then
\[
\gl=\frac k2,\qquad\theta=1 .
\]
\end{thm}

\begin{proof}
Let $w\in\Phi_X^{-1}(u_0)\cap\Theta$. By \eqref{eq:max-rank} and lower semicontinuity of the rank, $\dd \Phi_X$ has constant rank $k$ on a neighbourhood of $w$ in $\W_\RR$. It follows from the real analytic constant rank theorem that there are analytic diffeomorphisms $\alpha:(\RR^{\dim\W},0)\to(\W_\RR,w)$ and $\beta:(V_\RR,u_0)\to(\RR^M,0)$ with $\left(\beta\circ\Phi_X\circ\alpha\right)\left(x\right)=(x_1,\ldots,x_k,0,\ldots,0)$. Since $\beta$ is bi-Lipschitz near $u_0$, writing $x^{\prime}=(x_1,\ldots,x_k)$ we get, near $0$,
\[
c\,\|x^{\prime}\|^2\le\bigl\|\Phi_X(\alpha\left(x\right))-u_0\bigr\|^2\le C\,\|x^{\prime}\|^2,\qquad c,C>0 .
\]
Suppose first that $w\in\operatorname{int}\Theta$ then
$$\RLCT_w(I_{u_0};\varphi) \overset{(P2,P3)}{=} \RLCT_0\pr{(x_1,\ldots,x_k);1};$$
on the neighbourhood $\{\|x^{\prime}\|\le\epsilon\}\times\{\|x^{\prime\prime}\|\le\epsilon\}$, with $x^{\prime\prime}=(x_{k+1},\ldots,x_{\dim\W})$, polar coordinates in $x^{\prime}$ give
\[
\int\|x^{\prime}\|^{-s}\,\dd x^{\prime}\,\dd x^{\prime\prime}=c^{\prime}\int_0^{\epsilon}\rho^{\,k-1-s}\dd \rho ,
\]
which converges for $s<k$ and has a simple pole at $s=k$. Hence $\RLCT_w(I_{u_0};\varphi)=(k,1)$. If $w\in\partial\Theta$, the same computation gives $\RLCT_W(I_{u_0};\varphi)=(k,1)$ for a small neighbourhood $W$ of $w$, hence $$\RLCT_{\Theta,w}(I_{u_0};\varphi) \overset{(P4)}{\succeq} (k,1).$$ By hypothesis $\Phi_X^{-1}(u_0)\cap\operatorname{int}\Theta\neq\emptyset$, so the minimum in \eqref{eq:lambda-fibre} is $(k,1)$, that is $(2\gl,\theta)=(k,1)$.
\end{proof}

\begin{rem}\label{rem:generic}
By generic smoothness the set of $u\in\F(X)$ over which some parameter has rank smaller than $k$ is contained in a proper Zariski closed subset of $\F(X)$, so $(k/2,1)$ is the value of $(\gl,\theta)$ at a generic point of $\F(X)$. It is not necessarily the value at every smooth point: for the network $\RR\to\RR\to\RR$ with activation degree $r\ge2$ one has $\Phi_X(a,b)=ba^{r}x^{r}$, so $\F(X)=V$ is a line and every point is smooth, but at $u_0=0$ and $\Theta=[-1,1]^2$ the zeta function $\int|b|^{-s}|a|^{-rs}\dd a\dd b$ gives $(\gl,\theta)=\pr{\frac1{2r},1}$, while $\gl=\frac12$ at every $u_0\neq0$ realised in the interior. 
\end{rem}

For the converse we need coordinates in which $\Phi_X$ does not depend on the symmetry directions.

\begin{lem}\label{lem:slice}
Let $w\in\W^\circ_\RR$ and $u_0 = \Phi_X\left(w\right)$. There are real analytic coordinates $\left(a,z\right)\in\RR^{g}\times\RR^{d}$ on a neighbourhood of $w$ in $\W_\RR$, centred at $w$, such that
\[
\Phi_X\left(a,z\right)=\Psi\left(z\right),\qquad \Psi:(\RR^d,0)\to(V_\RR,u_0)\ \text{real analytic},\qquad \rk \dd \Psi_0=\rk \dd \Phi_{X,w};
\]
in particular $w$ is good if and only if $\dd \Psi_0$ is injective.
\end{lem}

\begin{proof}
In view of \cite[Theorem 1.16]{Gra26b}, the quotient map
\[
\pi_X:\W^\circ\longrightarrow\WW^\circ(X)
\]
is an étale $\G(X)$-torsor, hence a smooth morphism of relative dimension $g$, moreover $\WW^\circ(X)$ is smooth at $\pi_X\left(w\right)$ since $\W$ is smooth at $w$ and $\dim\WW^\circ(X)=\dim\W-\dim\G(X)=d$. Everything is defined over $\RR$, so near $w$ the map induced on real points is a real analytic submersion between manifolds of dimensions $\dim\W$ and $d$. By the submersion theorem there exist analytic coordinates $\left(a,z\right)$ centred at $w$, with $z$ centred at $\pi_X\left(w\right)$, in which $\pi_X\left(a,z\right)=z$. Since $\Phi_X=\phi_X\circ\pi_X$ \cite[Theorem 1.12]{Gra26b}, we get $\Phi_X\left(a,z\right)=\Psi\left(z\right)$, where $\Psi$ is $\phi_X$ read in the chart $z$. As $\dd\pi_{X,w}$ is surjective then $\rk \dd \Phi_{X,w}=\rk \dd \Psi_0$ and $w$ is good, that is $\rk \dd \Phi_{X,w}=d$, if and only if $\dd \Psi_0$ is injective.
\end{proof}

\begin{thm}\label{thm:B}
Let $\Theta$ be an admissible parameter domain, $\varphi$ an amplitude, and let
$w\in\W^\circ_\RR\cap\operatorname{int}\Theta$, with $u_0=\Phi_X(w)$ and
$r=\rk\dd\Phi_{X,w}$. Then the local learning coefficient at $w$ is at most $\frac{d+r}{4}$, and therefore
\[
\gl\le\frac{d+r}{4} .
\]
In particular, if $X$ is generically finitely identifiable and $w$ is not good, then $\gl<\frac d2=\frac k2$.
\end{thm}

\begin{proof}
Take the coordinates of Lemma \ref{lem:slice} and a product neighbourhood $Q\times B$ of $w$, with $Q\su\RR^g$ and $B\su\RR^d$ closed balls, contained in $\operatorname{int}\Theta$, and so small that $\|\Psi-u_0\|\le1$ on $B$, so that both integrals below are non-decreasing in $s$. 
If $\kappa:(\RR^g\times\RR^d,0)\to(\W_\RR,w)$ denotes the parametrisation given by the coordinates of Lemma \ref{lem:slice}, then the pulled back amplitude
\[
\varrho\left(a,z\right) = \varphi\pr{\kappa\left(a,z\right)}\,\bigl|\det \dd{\kappa_{\left(a,z\right)}}\bigr|,
\]
is continuous and bounded between two positive constants on $Q\times B$, so in view of Fubini's Theorem for real $s$
\[
\int_{Q\times B}\Bigl(\sum_j\pr{\Psi_j\left(z\right)-u_{0,j}}^2\Bigr)^{-s/2}\varrho\left(a,z\right) \dd a \dd z\ge c\int_{B}\Bigl(\sum_j\pr{\Psi_j\left(z\right)-u_{0,j}}^2\Bigr)^{-s/2}\dd z ,\qquad c>0 .
\]
By (P5) applied in $\RR^d$, and since an integral of a non-negative function diverges beyond its first pole, the right hand side diverges for $s\ge\frac{d+r}{2}$, hence so does the left hand side, and the first entry of $\RLCT_w(I_{u_0};\varphi)$ is at most $\frac{d+r}2$ and in view of \eqref{eq:lambda-fibre} we have $2\gl\le\frac{d+r}{2}$. If $k=d$ and $w$ is not good then $r<d$ and $\frac{d+r}{4}<\frac d2$.
\end{proof}

\begin{cor}\label{cor:dichotomy}
Let $X$ be generically finitely identifiable, let $\Theta$ be admissible, let
\[
u_0\in\Phi_X(\operatorname{int}\Theta)\setminus\Da(X)
\]
and let $\varphi$ be an amplitude. Then $(\gl,\theta)\preceq\pr{\frac d2,1}$. Moreover:
\begin{enumerate}
\item if every real parameter in $\Phi_X^{-1}(u_0)\cap\Theta$ is good, then $(\gl,\theta)=\pr{\frac d2,1}$;
\item if some real parameter in $\Phi_X^{-1}(u_0)\cap\operatorname{int}\Theta$ is not good, then $\gl<\frac d2$.
\end{enumerate}
\end{cor}

\begin{proof}
Since $u_0\notin\Da(X)$ the fibre lies in $\W^\circ$, and $k=d$. Let
$w\in\Phi_X^{-1}(u_0)\cap\operatorname{int}\Theta$. If $w$ is good, its local
threshold is $(d,1)$ by the proof of Theorem \ref{thm:A}; otherwise its first
entry is smaller than $d$ by Theorem \ref{thm:B}. In both cases
\eqref{eq:lambda-fibre} gives $(2\gl,\theta)\preceq(d,1)$. Items (1) and (2)
are Theorems \ref{thm:A} and \ref{thm:B}.
\end{proof}

\begin{rem}\label{rem:boundary}
The real fibre is a union of orbits of $\G(X)\left(\RR\right)$, which contains the unbounded rescaling tori, so it always meets $\partial\Theta$: this is why the hypothesis in Corollary \ref{cor:dichotomy}(1) is required on the whole of $\Phi_X^{-1}(u_0)\cap\Theta$. If the non-good real parameters over $u_0$ lie only on $\partial\Theta$, the value of $(\gl,\theta)$ depends on the shape of $\partial\Theta$ there, cf.~\cite[Example 2.8]{Lin17}. Goodness is $\G(X)$-invariant \cite[Proposition 1.22]{Gra26b}, so the hypothesis of (1) holds for every $\Theta$ as soon as one parameter is good in each $\G(X)$-orbit of the fibre that meets $\W_\RR$.
\end{rem}

\subsection{Exceptional singularities}
\label{Sec:exceptional}

\begin{cor}\label{cor:exceptional}
Let $X$ be generically finitely identifiable, let $u_0\in\E(X)$, and let $\Theta$ and $\varphi$ be as in Corollary \ref{cor:dichotomy}, with $u_0\in\Phi_X(\operatorname{int}\Theta)$.
\begin{enumerate}
\item If $u_0$ is not a real cusp in $\Theta$, then $(\gl,\theta)=\pr{\frac d2,1}$, the value at a generic point of $\F(X)$.
\item If $u_0$ is a real cusp witnessed by a parameter in $\operatorname{int}\Theta$, then $\gl<\frac d2$.
\end{enumerate}
In particular, by \cite[Theorem 2.7]{Gra26b}, an exceptional singularity can lower the learning coefficient only by being a cusp: exceptional nodes and exceptional points of non-properness which are not real cusps are invisible to $(\gl,\theta)$.
\end{cor}

\begin{proof}
Points of $\E(X)$ are outside $\Da(X)$, so this is Corollary \ref{cor:dichotomy} and Remark \ref{rem:generic}.
\end{proof}

The reason is that, in view of \eqref{eq:lambda-fibre}, $(\gl,\theta)$ is a minimum of invariants of the germs of $\Phi_X$ at the points of the fibre and, in view of the proof of Theorem \ref{thm:A}, at a good parameter that germ is, up to analytic changes of coordinates, given by a linear projection $(x^{\prime},x^{\prime\prime})\mapsto(x^{\prime},0)$. A node arises from the way several such branches meet in $V$, and a point of non-properness from branches whose parameters escape to infinity, which a compact $\Theta$ excludes; neither is seen by a single germ and only a rank drop is.

\begin{rem}\label{rem:block}
For the architectures of \cite[Remark 4.5]{Gra26b} the three parameters constructed in \cite[Construction 4.2]{Gra26b} are good, since each tangent image has dimension $\dim\WW(X)$; $w$ and $w^{\prime}$ defined there are real, so the local computation in the proof of Theorem \ref{thm:A} and (P1) give $\gl\le\frac d2$ at the corresponding exceptional singularity, for every $\Theta$ containing $w$ in its interior. Deciding equality requires the whole real fibre, which is not determined there.

By contrast, at the exceptional singularity of \cite[Section 3]{Gra26b} the whole fibre is known: it consists of three $\G(X)$-orbits of parameters, all of them good (see \cite[Section 3 and Remark 3.2]{Gra26b}), so Corollary \ref{cor:exceptional}(1) gives $(\gl,\theta)=(4,1)$ (since here $d=8$) there for every admissible $\Theta$ meeting the fibre in its interior.
\end{rem}

In two regimes the hypothesis of Corollary \ref{cor:exceptional}(1) is automatic, because every full parameter is good.

\begin{cor}\label{cor:regimes}
Assume that one of the following holds:
\begin{enumerate}
\item $X$ is shallow, with hidden layer $\N$, activation degree $r$ and $\dim X_{\mathrm{in}}\ge\dim\N\ge\dim X_{\mathrm{out}}$, and either $r\ge3$, or $r=2$ and $\dim\N=\dim X_{\mathrm{out}}$;
\item all activation degrees are equal to $r>9\pr{D(D-1)}^{L-1}$, where $D$ is the largest hidden width.
\end{enumerate}
Then $(\gl,\theta)=\pr{\frac d2,1}$ for every admissible $\Theta$, every $u_0\in\Phi_X(\operatorname{int}\Theta)\setminus\Da(X)$ and every amplitude. In particular every exceptional singularity realised in $\operatorname{int}\Theta$ is invisible to $(\gl,\theta)$.
\end{cor}

\begin{proof}
In case (1) every full parameter is good by \cite[Proposition 5.2]{Gra26b}. In case (2), a full parameter has no inactive and no redundant neuron, so it is not a subnetwork in the sense of \cite[Definition 3]{KMSSW26}, and is therefore a regular point of $\Phi_X$ by \cite[Theorem 5.1]{KMSSW26}, that is $\rk \dd \Phi_{X,w}=k$; moreover $k=d$ by \cite[(4)]{KMSSW26}. In both cases $X$ is generically finitely identifiable \cite[Corollary 1.24]{Gra26b}, and Corollary \ref{cor:dichotomy}(1) applies.
\end{proof}

By \cite[Remark 5.5]{Gra26b}, in case (1) the only possible exceptional singularities are points of non-properness, and in case (2) none is expected for networks in normal form \cite[Conjecture 6.7]{Gra26b}.




\appendix
\section{The hypotheses of Watanabe's theorem}\label{App:hypotheses}

We check the hypotheses of \cite[Chapter 6]{Wat09} and \cite[Chapter 6]{Wat18}
for the model \eqref{eq:network-density}, with $\Theta$ admissible and
$\varphi$ an amplitude. Throughout,
$e_1,\ldots,e_M$ is an orthonormal basis of $V_\RR$ and we write
\[
\delta(w)=\Phi_X(w)-u_0=\sum_{j=1}^M\delta_j(w)\,e_j ,
\]
so that each $\delta_j$ is a real polynomial on $\W_\RR$.

\subsection{Parameter set and prior.}
The set $\Theta=\{g_1\ge0,\ldots,g_l\ge0\}$ is compact and defined by finitely
many polynomial inequalities. The prior $\varphi$ is smooth on a neighbourhood
of $\Theta$ and positive on $\Theta$; since $\Theta$ has non-empty interior we have that $\int_\Theta\varphi>0$. Normalising $\varphi$ to a probability density changes
$F_n$ by an additive constant and does not affect $(\gl,\theta)$.

\subsection{Realisability and optimal parameters.}
The true distribution is $$q(\,\bullet\mid u_0)=p(\,\bullet\mid w_0),$$ so the
model is realisable. By \eqref{eq:network-KL}, proved below,
$K_{u_0}(w)\ge0$ with equality if and only if $\Phi_X(w)=u_0$. Hence the set
of optimal parameters is
\[
W_0=\{w\in\Theta:K_{u_0}(w)=0\}=\Phi_X^{-1}(u_0)\cap\Theta ,
\]
which is non-empty and compact.

\subsection{The log density ratio.}
Since $\log q(y\mid u)=-\frac M2\log(2\pi\sigma^2)-\frac{\|y-u\|^2}{2\sigma^2}$,
for every $w_0'\in W_0$ we have $p(\,\bullet\mid w_0')=q(\,\bullet\mid u_0)$, and
\[
f(y,w)=\log\frac{p(y\mid w_0')}{p(y\mid w)}
=\frac{\|y-\Phi_X(w)\|^2-\|y-u_0\|^2}{2\sigma^2},
\]
which does not depend on the choice of $w_0'\in W_0$. Let
$Y=u_0+\sigma\varepsilon$ with $\varepsilon\sim\mathcal N(0,\Id_M)$, so that $Y$
follows the true distribution, and let $\mathbb E$ denote expectation with
respect to it. Since
$\|\sigma\varepsilon-\delta\|^2-\|\sigma\varepsilon\|^2=-2\sigma\langle\varepsilon,\delta\rangle+\|\delta\|^2$,
\begin{equation}\label{eq:log-ratio}
f(Y,w)=-\frac1\sigma\sum_{j=1}^M\delta_j(w)\,\varepsilon_j+\frac1{2\sigma^2}\sum_{j=1}^M\delta_j(w)^2 .
\end{equation}
This is a polynomial in $w$ whose coefficients are constants or the random
variables $\varepsilon_j$, which lie in $L^s$ for every $1\le s<\infty$.
Hence $w\mapsto f(\,\bullet\,,w)$ extends to an entire map
$\CC^{\dim\W}\to L^s$, and in particular it is an $L^s$-valued real analytic
function on a neighbourhood of $\Theta$. Moreover, if
$m_\Theta=\max_\Theta\|\delta\|$, which is finite because $\Theta$ is compact,
the Cauchy-Schwarz inequality finally gives
\[
\sup_{w\in\Theta}|f(Y,w)|\le\frac{m_\Theta}{\sigma}\,\|\varepsilon\|+\frac{m_\Theta^2}{2\sigma^2}\in L^s .
\]

\subsection{The Kullback-Leibler divergence.}
Since $\mathbb E[\varepsilon_j]=0$ and $\mathbb E[\varepsilon_i\varepsilon_j]=\delta_{ij}$,
\[
\mathbb E\langle\varepsilon,\delta\rangle=\sum_j\delta_j\,\mathbb E[\varepsilon_j]=0,\qquad
\mathbb E\langle\varepsilon,\delta\rangle^2=\sum_{i,j}\delta_i\delta_j\,\mathbb E[\varepsilon_i\varepsilon_j]=\|\delta\|^2 .
\]
Taking expectations in \eqref{eq:log-ratio},
\[
K_{u_0}(w)=\mathbb E[f(Y,w)]=\frac{\|\delta(w)\|^2}{2\sigma^2},
\]
which is \eqref{eq:network-KL}.

\subsection{Relatively finite variance.}
Squaring \eqref{eq:log-ratio},
\[
f(Y,w)^2=\frac{\langle\varepsilon,\delta\rangle^2}{\sigma^2}
-\frac{\langle\varepsilon,\delta\rangle\,\|\delta\|^2}{\sigma^3}
+\frac{\|\delta\|^4}{4\sigma^4},
\]
and taking expectations as above,
\[
\mathbb E\bigl[f^2\bigr]=\frac{\|\delta\|^2}{\sigma^2}+\frac{\|\delta\|^4}{4\sigma^4}
=K_{u_0}(w)\pr{2+K_{u_0}(w)} .
\]
On $\Theta$ we have $K_{u_0}\le m_\Theta^2/2\sigma^2$, hence
\[
\mathbb E\bigl[f^2\bigr]\le\Bigl(2+\frac{m_\Theta^2}{2\sigma^2}\Bigr)K_{u_0}(w)
=\Bigl(2+\frac{m_\Theta^2}{2\sigma^2}\Bigr)\mathbb E[f] .
\]
Thus $\mathbb E[f]\ge c_0\,\mathbb E[f^2]$ for every $w\in\Theta$ and every
optimal parameter, with $$c_0=\pr{2+m_\Theta^2/2\sigma^2}^{-1}>0$$ and this gives
relatively finite variance \cite[Definition 7]{Wat18}.
\newpage

\bibliographystyle{amsplain}
\bibliography{References}

\end{document}